\documentclass[12pt]{amsart}
\usepackage{amssymb,amscd}

\newcommand{\A}{{\alpha}}
\newcommand{\D}{{\delta}}
\newcommand{\s}{{\sigma}}
\newcommand{\B}{{\beta}}
\newcommand{\cH}{{\mathcal H}}

\newcommand{\Z}{{\mathbb Z}}    
\newcommand{\R}{{\mathbb R}}    
\newcommand{{\C}}{{\mathbb C}}    

\def\g{{\gamma}}
\def\G{{\Gamma}}
\def\RA{{\rightarrow}}

\renewcommand{\section}[1]{\refstepcounter{section}\par\bigskip
\noindent\begin{center}{\normalsize\bf
\thesection .\hspace{3mm}{#1}}\end{center}\medskip \nopagebreak}

\def\theckbibliography#1{\par\bigskip
\begin{center}
{\normalsize \bf References}
\end{center}
\par
\noindent\list
 {[\arabic{enumi}]}{\settowidth\labelwidth{[#1]}\leftmargin\labelwidth
 \advance\leftmargin\labelsep
 \usecounter{enumi}}

 \sloppy\clubpenalty4000\widowpenalty4000
 \sfcode`\.=1000\relax}

\newtheorem{theorem}{Theorem}[section]

\newtheorem{corollary}[theorem]{Corollary}

\newtheorem{lemma}[theorem]{Lemma}

\newtheorem{proposition}[theorem]{Proposition}

\newenvironment{example}{\refstepcounter{theorem}\medskip\noindent{\bf
Example
\thetheorem}\hspace{1mm}}{\medskip}

\newenvironment{remark}{\refstepcounter{theorem}\medskip\noindent{\bf
Remark
\thetheorem}\hspace{1mm}}{\medskip}

\def\codim{{\mathrm{codim}}}
\def\Spec{{\mathrm{Spec}}}

\begin{document}

\begin{center}
\bigskip

{\Large\bf On the Smooth Points of $T$-stable \\
Varieties in $G/B$ and the Peterson Map }

\bigskip
\bigskip
\bigskip{\sc James B.\ Carrell\footnote{The first author was partially
supported by the Natural
Sciences and Engineering Research Council of Canada}\bigskip
\\Jochen Kuttler\footnote{The second author was partially supported by
the SNF
(Schweizerischer Nationalfonds)}}
\end{center}

\vspace{0.74in}
\begin{center}{\bf Abstract}
\end{center}

\bigskip
\bigskip
{\tiny Let $G$ be a semi-simple algebraic group  over ${\mathbb C}$, $B$
a
Borel subgroup of $G$ and
$T$ a maximal torus in $B$. A beautiful unpublished result of Dale
Peterson says that if
$G$ is simply laced,
then every rationally smooth point of a Schubert variety $X$ in $G/B$
is  nonsingular in $X$.
The purpose of this paper is to generalize this result to arbitrary
$T$-stable subvarieties
of $G/B$,  the only restriction being that  $G$ contains no $G_2$
factors.
A key idea in Peterson's proof is to deform the tangent space $T_y(X)$
to $X$ at a
nonsingular $T$-fixed point $y$
along the orbit of
$y$  under a root subgroup in $B$, which is open in a $T$-invariant
curve $C$ (we say
a $T$-curve) in $X$. In more
generality, if a
$T$-variety $X$ is  nonsingular along the open $T$-orbit in a $T$-curve
$C$ and $x\in C^T$,
we may consider the limit
$\tau_C(X,x)$ of the tangent spaces $T_z(X)$ as $z$ approaches $x$ along
$C$.
We call $\tau_C(X,x)$ the
{\it Peterson translate} of $X$ at $x$ along $C$. Peterson
showed that  a Schubert variety $X$ in $G/B$, where $G$ is semi-simple,
is nonsingular at
$x\in X^T$ as long as all
$\tau_C(X,x)$ coincide for all such {\em good} $T$-curves.. Our first
result
generalizes this theorem to any irreducible $T$-variety $X$,
provided the fixed point $x$ is attractive under much weaker hypotheses.
We then prove that if
$X$ is a $T$-variety in $G/B$  where $G$ contains no
$G_2$ factors, then every Peterson translate $\tau_C(X,x)$ is contained
in the
linear span $\Theta_x(X)$ of the
reduced tangent cone to $X$ at $x$. (This fails when $G=G_2$.) Combining
these
two results leads  to our
characterization of the smooth $T$-fixed points of such
$X$. In particular, we show that if $G$ is simply laced, then $X$ is
nonsingular
at a $T$-fixed point $x$  if
and only if it is rationally smooth at $x$ and $x$ lies on at least two
good $T$-curves
Peterson's $ADE$ result is an immediate consequence of this. In
addition, we
obtain a uniform description of
the nonsingular
$T$-fixed points of a Schubert variety in $G/B$ modulo the $G_2$
restriction.
In particular, a Schubert variety $X$
in such a $G/B$ is nonsingular if and only if all the reduced tangent
cones of $X$ are linear.
Finally, we also obtain versions of the above results for all algebraic
homogeneous spaces $G/P$ modulo
the $G_2$ restrictions.}

\newpage

\section{Introduction}\label{intro}

Let $G$ be a semi-simple algebraic group over $k={\mathbb C}$. Fix a
Borel subgroup $B$ of
$G$ and a maximal torus $T\subset B$.
The purpose of this paper is to investigate the singular locus  of a
$T$-stable subvariety $X$ of the flag variety $G/B$. More precisely, we
would
like to describe the set of nonsingular $T$-fixed points of
$X$. This problem originates with the question of determining the
connection
between the singular loci of a Schubert variety $X\subset G/B$ (i.e. the
closure of a $B$-orbit) in the sense of  rational smoothness  (cf.
\cite{kl1,kl2}) and the sense of algebraic geometry. It was shown, for
example, in
\cite{deod} that when
$G$ is of type $A$, i.e. $G/B$ is the variety of complete flags in
$k^n$, then
the two singular loci are the same, in particular every rationally
smooth point of $X$ is
nonsingular. More recently, Dale Peterson (unpublished) extended this to
Schubert varieties in
$G/B$ in the full $ADE$ setting (see below). His method is to study how
tangent
spaces of a Schubert variety $X$ behave when deformed along
$T$-invariant curves in $X$ containing a nonsingular point of $X$.

Before describing our results, we need to fix some notation. Recall that
the
$T$-fixed point set $G/B^T$ is in a one to one correspondence with the
Weyl group
$W$ of $(G,T)$ via $w\mapsto wB$. Hence we may simply denote
$wB\in G/B^T$ by $w$. The Schubert variety $X(w)$ associated to
$w\in W$ is by definition the Zariski closure in $G/B$ of the $B$-orbit
$Bw$.
Recall that $B$ defines a Coxeter system for $W$. Let
$\le$ denote the
associated partial order on $W$, the so called Bruhat-Chevalley order.
This
Coxeter system has two fundamental properties. Firstly, $x\leq y$ if and
only
if $X(x)\subset X(y)$. Hence
$X(w)^T=\{x\leq w\}$. Note that we will usually use $[x,w]$ to denote
$\{x\leq w\}$.
Secondly, if $\ell(w)$ denotes the length of $w\in W$, then
$\ell(w)=\dim X(w)$.

For simplicity, let $X$ denote $X(w)$. The set $E(X,x)$ of $T$-invariant
curves, or
briefly, $T$-curves, in $X$  containing the $T$-fixed point $x$ turns
out to be of basic
importance in determining
the singular locus of $X$  (cf. \cite{cp,spsv}). Let $\Phi\subset X(T)$
be the root system
of $(G,T)$, and recall that to each $\A\in \Phi$, there is a one
dimensional unipotent subgroup
$U_\A$ of $G$ called the {\em root subgroup}
associated  to $\A$. Recall that the positive roots $\Phi^+$
can be described as those such that $U_\A \subset B$. Then any $C \in
E(X,x)$  has the form $\overline{U_\A x}$ for some $\A$. Moreover,
$C^T=\{x,y\}$, where  $y=r_\A x$, $r_\A$ denoting the reflection
corresponding
to $\A$. When $y>x$, then $\A < 0$ and we can write $C = \overline{U_\B
y}$ with $\B = -\A >
0$, so one can translate the Zariski tangent space $T_y(X)$ to $X$ at
$y$ along
$C\backslash \{x\}$ via $U_{\B}$ leaving $X$ invariant. Taking the limit
gives a $T$-stable subspace
$\tau_C(X,x)$ of $T_x(X)$ of dimension $\dim T_y(X)$. The key result is

\bigskip\noindent
{\bf Peterson's Theorem} {\em Suppose that $X=X(w)$ is nonsingular at
all $y\in W$ such that $x<y\leq w$ and that all $\tau_C(X,x)$ coincide
when
$C\in E(X,x)$ has the property that $X$ is nonsingular on
$C\backslash \{x\}$. Then
$X$ is nonsingular at $x$.}

\bigskip
The idea of Peterson's proof is to show that if all the $\tau_C(X,x)$
coincide, then the fibre over
$x$ in the Nash blow up of $X$ at $x$ contains no $T$-curves. Since
Schubert
varieties are normal, it
follows from Zariski's Connectedness Theorem and Lemma
\ref{TC} that this fibre consists of a single point. This implies, by a
result of Nobile \cite{no}, that
$X$ is nonsingular at $x$.  Using this, Peterson was able to show
\begin{theorem}\label{pade} If $G$ is of type $ADE$, then every
rationally smooth point  of a Schubert variety  $X(w)$ in $G/B$ is
nonsingular.
\end{theorem}

Combining this result with the characterizations of rationally smooth
Schubert varieties given in \cite{cp}, we get several lovely
descriptions of the nonsingular Schubert
varieties in $G/B$ for the simply laced setting.

\begin{theorem} {\rm (cf Theorem A of \cite{cp})} Let $G$ be
semi-simple. Then a Schubert
variety
$X(w)$ in $G/B$  is  rationally smooth if and only if
any of the following equivalent conditions hold:

\begin{itemize}

\item[\rm{(1)}] the Poincar\'e polynomial of $X(w)$
$$P(X(w),t)=\sum b_{i}(X(w))t^i=\sum_{x\le
w}t^{2\ell(x)}$$
is symmetric;

\medskip
\item[\rm{(2)}] the order $\le$ on $[e,w]$ is rank symmetric;

\medskip
\item[\rm{(3)}] for each $x\le w$, $|E(X(w),x)|=\ell(w)$; and

\medskip
\item[\rm{(4)}] the average $a(w)$ of the length function on $[e,w]$ is
$\frac{1}{2}\ell(w)$. that is,
$$a(w)=\frac{1}{|[e,w]|}\sum_{x\le w} \ell(x)=\frac{1}{2}\ell(w).$$

\end{itemize}
\end{theorem}

We therefore obtain
\begin{corollary} \label{RSGmodP} If $G$ is simply laced, then a
Schubert variety $X(w)$ in $G/B$ is
nonsingular  if and only if any of the equivalent conditions {\rm
(1)-(4)} hold.
\end{corollary}

A  corresponding $G/P$ version will be stated in \S
\ref{G/P}.

We now describe some  generalizations of these results for arbitrary
irreducible  $T$-stable subvarieties $X$ of $G/B$ proved in this
paper.  Put
$$TE(X,x)=\sum_{C\in E(X,x)} T_x(C).$$
If $C=\overline{U_{\A}x}$, then $T_x(C)$ is a $T$-stable line of weight
$\A$, so  $TE(X,x)$
is a $T$-submodule of $T_x(X)$ such that
$\dim TE(X,x)=|E(X,x)|$
(cf. \cite{cp}). In particular $\dim T_x(X)\geq
|E(X,x)|$. We will call
$C\in E(X,x)$ {\em  good} if $X$ is nonsingular along the open $T$-orbit
in $C$, or, equivalently, if
$C$ is not contained in the singular locus of $X$.  Our generalization
of Peterson's  Theorem goes as
follows:
\begin{theorem}\label{TE} Suppose $\dim X\geq 2$ and $x\in X^T$. Then a
necessary
and sufficient condition
that  $X$ be nonsingular at $x$ is that there exist at least two
distinct good $T$-curves
$C,D\in E(X,x)$ such that
\begin{equation}
\tau_C(X,x)=\tau_D(X,x)=TE(X,x).
\end{equation}
If $X$ is Cohen-Macaulay  at $x$, then $X$ is nonsingular at $x$ if and
only if there exists
at least one good $C\in E(X,x)$ such that $\tau_C(X,x)=TE(X,x)$.
\end{theorem}

The proof only uses the Zariski-Nagata Theorem, and hence is
completely algebraic. In particular, it works over any
algebraically closed field. (This improves the proof given in
\cite{kut}.) If $X$ is a Schubert variety, it is not hard to show
that Theorem \ref{TE} implies Peterson's Theorem, giving the
first of four several proofs. The Cohen-Macaulay  statement is
proved in Proposition \ref{CM}.

In order to generalize Peterson's $ADE$ Theorem, we need to understand
where the
$\tau_C(X,x)$ are situated in $T_x(X)$. For
this, let $\Theta_x(X)$ denote the linear span of the reduced tangent
cone
of $X$ at $x$. If $C\in E(X,x)$ has the form $C=\overline{U_\A x}$, we
will call $C$
{\it long} or {\it short} according to whether $\A$ is long or short. If
$G$ is simply laced, then, by
convention, all
$T$-curves will be called short. Clearly,
$$TE(X,x)\subset \Theta_x(X)\subset T_x(X).$$
The next result is one of our key observations.

\begin{theorem}\label{tausubsetT} Assume $G$ has
no $G_2$ factors. Then, if $C\in E(X,x)$ is good,
$$\tau_C(X,x)\subset \Theta_x(X).$$
Moreover, if $C$ is short, then $\tau_C(X,x)\subset TE(X,x)$. In
particular, if $G$ is simply laced and $C$ is good, then
$\tau_C(X,x)\subset TE(X,x)$
\end{theorem}
We give an example in \S \ref{tancone} which shows that the $G_2$
hypothesis is necessary.
Peterson's $ADE$ Theorem is now a simple consequence. Indeed, assuming
$X=X(w)$, it suffices to
suppose $x$ is a rationally smooth $T$-fixed point such that $X$ is
smooth at every $y$ with $x<y\leq
w$. Since the singular locus of a Schubert variety has codimension at
least two,
$\ell(x)\leq \ell(w)-2$. Hence $x$ lies on at least two good $T$-curves
(cf. Proposition 2.3). The proof now
follows from Theorems \ref{TE} and \ref{tausubsetT}, since if $X$ is
rationally smooth at $x$, then
$|E(X,x)|=\dim X$.

This result now gives us a complete description of the smooth points of
Schubert varieties in $G/B$ as long as $G$
contains no $G_2$ factors.
\begin{corollary}\label{SPSV} Suppose  $G$ contains no $G_2$ factor.
Then the Schubert variety $X=X(w)$ is smooth
at $x<w$ if and only if $\dim \Theta_y(X)=\dim X$ for all $y\in [x,w]$.
In other words, $X$ is smooth at
$x$ if and only if the reduced tangent cones of $X$ at all $y\in [x,w]$
are linear.
In particular, $X$ is smooth if and only if all its reduced tangent
cones are linear.
\end{corollary}

The proof is essentially the same as that of the above proof of the
$ADE$ Theorem.

A natural question is whether Peterson's $ADE$ Theorem holds for
arbitrary
$T$-varieties in $G/B$ if $G$ is simply
laced. It turns out that the answer is in general no, but we do have the
following:
\begin{corollary} If $G$ is simply laced, $X$ is rationally smooth at
$x$
and $\dim X\geq 2$, then $X$ is smooth at
$x$ if and only if $E(X,x)$ contains at least two good $T$-curves.
\end{corollary}

Indeed, if $X$ is rationally smooth at $x$, then by a recent result of
Brion \cite{bri},
$|E(X,x)|=\dim X$. Hence
the corollary follows immediately from Theorem \ref{TE}.

If $X$ is a  Schubert variety and $G$
is simply  laced, then we know from \cite{cp,car} that
$\Theta_x(X)=TE(X,x)$. In
fact, if $G$ is simply laced, this turns out to be true for all
irreducible
$T$-subvarieties of $G/B$.
\begin{theorem} \label{sltc} Assume $G$ is simply laced, and $x\in X^T$.
Then every $T$-line
in the reduced
tangent cone to $X$ at $x$ has the form $T_x(C)$ for some $C\in E(X,x)$.
That is,
$$\Theta_x(X)=TE(X,x).$$
\end{theorem}

We now briefly describe the rest of the paper.
First of all,  in \S \ref{PM}, we define the Peterson map in a general
setting and derive
its basic properties. In particular, if $X$ is a Schubert variety,  we
show there is a
remarkable explicit formula for
$\tau_C(X,x)$ for any (not necessarily good) $C\in E(X,x)$. In \S
\ref{FL}, we prove a
fundamental lemma showing that the Peterson map for a good
$T$-curve $C$ is completely determined by its behavior on the
$T$-surfaces in $X$ containing $C$. This gives us Theorem
\ref{tausubsetT} and,
in addition, allows us to deduce that certain weights outside $TE(X,x)$
may occur in
$\Theta_x(X)$. However, the fact that there is no general description of
$\Theta_x(X)$ makes it desirable to find a subspace containing
$\tau_C(X,x)$, assuming $C$ is good, admitting an explicit description.
We describe
such a subspace in \S \ref{PTSV}. In
the next section, we mention an algorithm for finding the singular locus
of a Schubert
variety in $G/B$. In \S \ref{G/P}, we prove a lemma which extends our
results to any $G/P$,
with the suitable restrictions on $G$, and in the last  section we
mention some open problems.

A remark about the field is in order. Although we are assuming
$k={\mathbb C}$, we believe our
arguments are valid over any algebraically closed field. This goes hand
in hand
with the fact proved in \cite{po} that the singular locus of a Schubert
variety is independent of the
field of definition.

\bigskip
\noindent{\bf Acknowledgement} The authors would like to thank Dale
Peterson for
discussions about his results. We also thank Michel Brion for some
comments on rational smoothness.

The second author would like to thank the University of British
Columbia, Vancouver, and the University
of California, San Diego, for hospitality during the work on this paper.

\section{Preliminaries on $T$-varieties}\label{prelim}

Throughout this paper, $T$ will  denote an algebraic torus over $k=
\mathbb C$ with character group $X(T)$
and dual group $Y(T)$ of one parameter subgoups of $T$. $X$ will always
denote an irreducible
$T$-variety with finite non-empty fixed point set $X^T$ which is
locally linearizable in the following sense: every point $z\in X$ has a
connected affine $T$-stable neighborhood $X_z$  admitting a
$T$-equivariant embedding into an affine space $V$ with a linear
$T$-action. This is for example true for closed $T$-stable subsets of a
normal $T$-variety (\cite{sumi1},\cite{sumi2}). For any $T$-variety $X$
and $x \in X$ we choose once and for all such a neighborhood $X_x$.

We will denote the set of weights of a  $T$-module $V$ by $\Omega (V)$.
If
$x\in X^T$ and  all elements of
$\Omega (T_x(X))$ lie on one side of a
hyperplane in $X(T)\otimes \mathbb Q$, then $x$ is called  {\em
attractive}. It follows immediately from the definition, that if $x$ is
attractive,
there exists a one
parameter group $\lambda \in Y(T)$ such that $\langle \A, \lambda
\rangle >
0$, where $\langle \cdot, \cdot \rangle :
X(T) \times Y(T) \rightarrow \mathbb Z$, is the natural pairing.
Equivalently, $x$ is attractive if and only if $\lim_{t \rightarrow
0}\lambda(t)y
= x$ for all $y \in X_x$. It is well known that if $x \in X^T$ is
attractive, there is a closed $T$-equivariant immersion $X_x
\subset T_x(X)$. For example, every $T$-fixed point in  $G/B$ is
attractive.
If $x$ is attractive and  $L\subset T_x(X)$ is a $T$-stable line, we may
consider the restriction to $X_x$ of a $T$-equivariant linear projection
$T_x(X) \RA L$. Since $L$ is an affine line, this restriction gives rise
to a $T$-eigenvector $f \in k[X_x]$.  We say that $f \in k[X_x]$ {\em
corresponds} to $L$ if $f$ is so obtained.

Another fact, that we will use below, is
\begin{lemma} \label{AFF} Let $X$ be affine and  $x\in X^T$ attractive.
If $Y$ is any affine $T$-variety, then a $T$-equivariant morphism $f: X
\rightarrow Y$ is finite if and only if $f^{-1}(f(x))$ is a finite set.
\end{lemma}

As in the introduction, $E(X,x)$ will denote the set of
$T$-curves in $X$ containing the point $x\in X^T$. Also as above, a
$T$-curve $C$ is called \emph{good} if $C^o=C\backslash C^T\subset
X^*$, where $X^*$ is the set of nonsingular points in $X$. This just
means that $C\cap X^*$ is nonempty.
The following
lemma gives a very useful fact about
$E(X,x)$ (cf \cite{cp}).

\begin{lemma}\label{TC}
For any $x\in X^T$,
$$|E(X,x)|\geq \dim_x X.$$
That is, the number of $T$-curves in $X$ through $x$ is at least $\dim_x
X$.
\end{lemma}


If the number of $T$-curves in $X$ is finite, then there is a finite
graph $\G(X)$, called the {\em Bruhat graph} of the pair $(X,T)$,
which generalizes the Bruhat graph of the
Weyl group $\G(W)$ (see for example
\cite{cp}). The vertices of
$\G(X)$ are the $T$-fixed points, and two $x,y\in X^T$ are joined by an
edge if and only if there
exists a $T$-curve $C$ in $X$ such that $x,y\in C$. When $X$ is a
$T$-variety in $G/B$, where $T$ is
maximal in $B$, then $\G(X)$ is a subgraph of $\G(W)$. In
particular, if $x,y\in X^T$ are joined by an edge
in $\G(X)$, then there exists an $r\in R$ such that $x=ry$. If $X$ is a
Schubert
variety, then  $\G(X)$ is a full subgraph; any edge of $\G(W)$ joining
$x,y\in X^T$ is also an edge of
$\G(X)$. Notice also that the Chevalley-Bruhat order gives a natural
order on the vertices of
$\G(X)$. In the Schubert case, Lemma \ref{TC} implies
\begin{proposition} \label{DI} \rm{(Deodhar's Inequality} \cite{cp}) If
$x<w$,
then there
exist at least $\ell(w)-\ell(x)$ reflections $r\in W$ for which
$x<rx\leq w$.
\end{proposition}

\section{The Peterson Map}\label{PM}

Let $X$ be a $T$-variety, $x\in X^T$ a locally linearizable isolated
fixed
point, and assume $C\in E(X,x)$. In this  section, we will define and
study what we call the {\em Peterson map} $\tau_C(~,x)$. The version we
consider here is slightly more general than the
tangent space deformation considered by Peterson, which was only defined
in the case of Schubert varieties. In the next section, we will relate
these two deformations and give an explicit computation Peterson's
version.

Our Peterson map is defined on certain subspaces of the tangent space to
$X$ at an
arbitrary point $z\in C^o$ of a $T$-curve $C$ in $X$.
We may suppose the $T$-stable neighborhood $X_x \subset V$ is embedded
equivariantly
into a $T$-module $V$. Let $M\subset
T_z(X)$ be a
$k$-subspace stable under the isotropy group $S$ of $z$, and let
${\mathbf
M}^o = T M \subset T(X) \left |
_{C^o}\right . $ be the
$T$-stable vector-bundle over $C^0$, having fibre $M$ over $z$. Define
${\mathbf
M}$ to be the Zariski closure of ${\mathbf
M}^o$ in $T(X)$. Then, by definition, the Peterson map assigns
$\tau_C(M,x)={\mathbf  M}\cap T_x(X)$ to $M$. If
$M=T_z(X)$, then we will denote $\tau_C(M,x)$ by
$\tau_C(X,x)$. Clearly, if $X$ is nonsingular at $x$, then
$\tau_C(X,x)=T_x(X)$.

If $C$ is smooth, then an alternative description of $\tau_C(M,x)$ is as
follows.
By the properness of Grassmannians, the $T$-stable  vector bundle
${\mathbf
M}$ on $C\cap X_x$ extends to a vector bundle $\mathbf M$ on $C$ such
that
the
restriction
of ${\mathbf M}$ to $C^o$ is ${\mathbf M}^o$. Then $\tau_C(M,x)={\mathbf
M}_x$.

  The main properties of the Peterson map are given in the
next result. We assume the notation defined above is still in effect.
\begin{proposition} \label{PPP} Suppose $X$, $T$, $x\in X^T$ and $C\in
E(X,x)$ are as above, and let
$M$ be an $S$-stable subspace of $T_z(X)$ of dimension $m$. Then:

\begin{itemize}
\item[\rm{(1)}] $\tau_C(M,x)$ is a $T$-stable subspace of $T_x(X)$ of
dimension
$m$, and, moreover, $M$ and
$\tau_C(M,x)$ are isomorphic $S$-modules. If $Y \supset X$ is any
ambient
smooth $T$-variety, then as elements of the Grassmannian ${\mathcal
G}_m(Y)$ of
$m$-planes in $T(Y)$, $\tau_C(M,x) = \lim_{t\RA 0}{\mathrm
d}\lambda(t)M$,
where $\lambda$ is an arbitrary one parameter subgroup of $T$ such that
$\lim_{t\RA 0 }\lambda(t)z=x$.

\item[\rm{(2)}] If $M=M_1\oplus \dots \oplus M_t$ is the $S$-weight
decomposition of $M$, then $\tau_C(M,x)= \tau_C(M_1,x)\oplus \dots
\oplus
\tau_C(M_t,x)$ is the $S$-weight decomposition of $\tau_C(M,x)$.

\item[\rm{(3)}] If $N$ is any $T$-stable subspace of $\tau_C(X,x)$, then
there exists an $S$-stable subspace $M\subset T_z(X)$ such that
$\tau_C(M,x)=N$.
\end{itemize}
\end{proposition}

\begin{proof}
It follows from the definitions that $\tau_C(M,x)$ is $T$-stable. That
it
is a subspace of the same
dimension as $M$ also follows from the properness of the Grassmannian.
Moreover, it follows easily that the limit in the Grassmannian is found
by closing the bundle $T(X)$ over
$C^o$. Given all this, the first statement follows from the second, once
we have proved that
$\tau_C(M_i,x)$ and $M_i$ are isomorphic as $S$-modules. Now $S$ acts on
$M_i$
by a character $\A_i$, hence on $TM_i \subset T(X)$, as well. With
$TM_i$
being dense in $\overline{TM_i}$, it follows that $sv = \A_i(s)v$ for
all
$v \in \tau_C(M_i,x), s \in S$. It is now obvious that  $\tau_C(M,x)$
decomposes as stated. Thus we have 1) and 2).

For the last statement, it is enough, of course, to consider the case
where $N$ is a line having $T$-weight say $\A$. By 1) and 2), there is
an
$S$-subspace $M_1$ of $T_z(X)$, on which $S$ acts by the character $\A$
restricted to $S$, so that $\tau_C(M_1,x)$ contains $N$. Let $\lambda
\in
Y(T)$
be a regular one parameter subgroup of $T$ so that $\lim_{t\rightarrow
0}\lambda(t)z=x$. Thus, there is an induced surjective morphism $f:
{\mathbb
A}^1 \RA C_x \subset X_x$. Let $V$ be a $T$-module, into which $X_x$
embeds
equivariantly with $x = 0$. Then
$${\mathbf B} = f^*(T(X_x)) ={\mathbb A}^1\times_{C_x} T(X_x)$$
is a $T$-stable  subvariety of ${\mathbb A}^1\times V$,
which certainly contains ${\mathbf M}^\prime = \overline{f^*(TM_1)}$,
the
latter being a vector bundle over ${\mathbb A}^1$. This means in
particular
that ${\mathbf M}^\prime$ is a trivial bundle. Moreover, it is easy to
see,
that over $
{\mathbb G}_m \subset {\mathbb A}^1$, the map
$(s,w)\mapsto(s,{\mathrm{d}}\lambda (s)w)$ is a closed
$S$-equivariant immersion, $S$ acting trivially on the first factor
${\mathbb A}^1$.

Summarizing, there is a global section $\sigma$  of ${\mathbf
M}^\prime$
with
$\sigma(1) \in M_1$ and with $0 \not= \sigma(0) \in N$, which, over
${\mathbb G}_m$ has the form
$$\sigma(s)=(s, {\mathrm d}\lambda(s)(\sum_i s^i w_i)),$$
for suitable $w_i \in M_1 \subset V$. Let $V = \bigoplus_\B V_\B$ be
the decomposition into $T$-weightspaces, so that $w_i = \sum_\B
v_{i,\B}$
with $v_{i,\B} \in V_\B$. Now compare weights in the expansion
$$\sum_i{\mathrm d}\lambda(s)s^{i}w_i = \sum_{i,\B} s^{i + \langle \B ,
\lambda \rangle}v_{i,\B}.$$ Since $\sigma$ extends to zero and
$\s(0)\in N\setminus \{0\}$, the term on the right hand side of
(2) of degree zero occurs when $i = - \langle\A,\lambda \rangle$.
Furthermore, for every nonzero term on the right hand side,
$i+\langle \B, \lambda\rangle \geq 0$. Since $\lambda$ is
regular, it follows that $\tau_C(kw_i,x)=kv_{i,\A}$, where $i = -
\langle\A,\lambda\rangle$. This says $\tau_C(kw_i,x)=N$, and we
are done.
\end{proof}

\section{The Peterson Map for $G/B$}\label{G/B}

In this section, we will explicitly compute  $\tau_C(X,x)$ for a
Schubert variety $X=X(w)$ in $G/B$, where $C\in
E(X,x)$ is such that $C^T=\{x,y\}$ and $y>x$. In other words, we are
considering what happens when we pass from a
higher vertex of the Bruhat graph along an edge to a lower vertex. Then
$C$ can be expressed in the form
$C=\overline{U_{\A}y}$ with $\A>0$, hence  the additive group
$U_\A$ acts transitively on
$C\backslash \{x\}$. Thus any subspace $M\subset T_zX$,
$z\in C^o$, is the $U_\A$-translate of a unique subspace of
$T_yX$.  In addition, the vector bundle
${\mathbf M}$ introduced in the previous section is  defined  and
$U_\A$-equivariant on $C$. Therefore, $\tau_C$ can
be viewed as defined on subpaces of $T_y(X)$.  This is  the map
originally considered by Peterson.

Letting $S=\ker (\A)$, suppose $M\subset T_y(X)$ is $S$-stable (resp.
$T$-stable). Then
${\mathbf M}$ is also $S$-equivariant (resp. $T$-equivariant). Thus
$\tau_C(M,x)$ is an $S$-module (resp. $T$-module), and furthermore,
$\tau_C(M,x)$ is also $U_\A$-module. (This does not require that
$M$ be $S$-stable.)
Assuming $M$ is $S$-stable, any $S$-weight space $V$ of $M$ is a direct
sum of certain ${\mathbf g}_{\B+k\A}$, $\B$ is fixed,  $k\leq 0$ and
$y^{-1}(\A+k\B)<0$ (since
$\Omega (T_y(G/B))=y^{-1}(-\Phi^+)$).

Now suppose $M$ is an $S$-weight space in $T_y(X)$ of dimension $\ell$.
Then $M$
and $\tau_C(M,x)$ are isomorphic as
$S$-modules, but possibly different when viewed as $T$-modules. However,
the
$T$-weights of $\tau_C(M,x)$ are not hard to determine. Indeed, since
$M$
has only one $S$-weight, it follows that $\Omega (M)$ is contained in a
single
$\A$-string in $\Phi$. Now $T_y(G/B)$ is a ${\mathbf g}_{-\A}$-module,
although
$T_y(X)$ need not be one. In fact there exists a unique ${\mathbf
g}_{-\A}$-submodule $M^*$ of $T_y(G/B)$ such that
$M^*\cong M$ as $S$-modules. That $M^*$ exists is clear. If $M$ isn't
already a ${\mathbf g}_{-\A}$-module, then one way to describe $M^*$ is
as the
unique $U_{-\A}$-fixed point on the $T$-curve $\overline{U_{-\A}M}$
in the ordinary Grassmannian $G_{\ell}(T_y(G/B))$. Clearly $M^*$ is
determined by
the unique $\B\in \Phi$ lying on the $\A$-string containing $\Omega
(M)$
satisfying the condition that $y^{-1}(\B-\ell \A)\not\in \Phi^-$, but
$y^{-1}(\{\B, \B-\A, \dots ,\B-(\ell -1)\A \})\subset \Phi^-$. That is,
$$\Omega (M^*)=\{\B, \B-\A, \dots ,\B-(\ell -1)\A \}.$$
We will call
$\B$ the {\em leading weight} of $M^*$.
\begin{proposition}\label{PMG/B} Assuming $M$ is as above,
$\tau_C(M,x)={\mathrm{d}} \dot{r}_\A(M^*)$, where
${\mathrm{d}} \dot{r}_\A$ denotes the differential at $y$ of a
representative  ${\dot{r}}_\A\in N(T)$ of
$r_\A$. Consequently, $\Omega (\tau_C(M,x))=r_{\A}\Omega (M^*)$.
\end{proposition}

\begin{proof}
Since $M^*$ is a $U_{-\A}$-module, $\dot{r}_A(M^*)$ is
a $U_\A$-module contained in $T_x(G/B)$ which is isomorphic to $M$ as an
$S$-module. But this condition uniquely determines $\tau_C(M,x)$.
\end{proof}

\begin{remark} Of course we do not actually need to define $M^*$ to
determine
$\tau_C(M,x)$. We will use this formulation in \S \ref{SL}.
\end{remark}

Recall that if $X(w)$ is smooth at $x\leq w$, then
$$\Omega (T_x(X(w))=\Omega (TE(X,x))=\{\g\in \Phi\mid x^{-1}(\g)<0,~
r_\g x\leq w\}.$$
If $x<y\leq w$, then clearly $X(w)$ is smooth at $y$ also. Using the
Peterson map, we can now describe how to
obtain the $T$-weights of $T_x(X(w))$ by degenerating to $x$ along
edges of $\G(X(w))$.
Denoting $\Omega (T_x(X(w)))$ by $\Phi(x,w)$ as in
\cite{cp} and putting
$\Phi(x,w)^*=\Omega (T_x(X(w))^*)$, Proposition \ref{PMG/B} gives
\begin{corollary} Let $X(w)$ be smooth at two adjacent vertices of
the Bruhat graph $\G(X(w))$, say $x$ and $y=rx$, where
$y>x$ and $r\in R$. Then
$$\Phi(x,w)=r(\Phi(y,w)^*).$$
Consequently, if $x<z\leq w$ is another vertex of $\G(X(w))$ adjacent to
$x$,
then
$r(\Phi(y,w)^*)=t(\Phi(z,w)^*)$, where $x=t z$ with $t\in R$.
\end{corollary}
We will see later that if $X(w)$ is smooth at $y$ but not necessarily at
$x=ry<y$, then it is still true that
$r(\Phi(y,w)^*)\subset \Phi(x,w),$ provided the corresponding $T$-curve
$C\in
E(X,x)$  is short. Note: we are assuming (contrary to the common
practice) that if $\Phi$ is
simply laced, then all its elements are short. Thus all $T$-curves in
the
corresponding $G/B$ are by convention short.

If $C$ is long, the
situation is more complicated. This is illustrated in the
following example.

\begin{example}\label{B2} Suppose $G$ is of type $B_2$, and $w=r_\A r_\B
r_\A$,
where $\A$ the short simple root and $\B$ the long simple root.
In this example, we will compute the Peterson maps and use the result to
 determine the singular locus of $X(w)$, which is of course already well
known. Put $X=X(w)$ and
$\Omega (T_x(X))=\Omega (x)$.
If $x\leq w$  is a nonsingular point of $X$, then $\Omega
(x)=\Phi(x,w)$.
Clearly (for example, by Peterson's Theorem), $w$, $r_\A r_\B$ and $r_\B
r_\A$
are nonsingular points, and
one easily sees that
\begin{itemize}
\item [(1)] $\Omega (w)=\{\A,\A+\B,2\A+\B\}$;
\item [(2)] $\Omega (r_\A r_\B)=\{\A,2\A+\B,-(\A+\B)\}$;
\item [(3)] $\Omega (r_\B r_\A)=\{-\A,\B,\A+\B\}$;
\end{itemize}
It remains to test whether the points $r_\A$ and $r_\B$ are nonsingular.
Indeed,
since $\A$ is simple and $r_\A w<w$, $\dot{r}_\A X=X$. Moreover, if
$C=\overline{U_\A x}$, then $\Omega (\tau_C(X,x))=r_\A(\Omega (r_\A x))$
as long as
$r_\A x<
x\leq w$. Thus
$e$ is a nonsingular point if and only if $r_\A$ is. Let's first compute
$\Omega (\tau_C(X,r_\A)$ where
$C=\overline{U_\B r_\B r_\A}$.
It is clear that $\Omega (r_\B r_\A)$ is the set of weights of
of a ${\mathbf g}_{-\B}$-submodule of $T_{r_\B r_\A}(G/B)$, so
$$\Omega (\tau_C(X,r_\A))=r_\B (\Omega (r_\B
r_\A))=\{-(\A+\B),-\B,\A\}.$$
Next consider $\Omega (\tau_D(X,r_\A))$ where
$D=\overline{U_{2\A+\B} r_\A r_\B}.$
It is again clear that $\Omega (r_\A r_\B)$ is the set of weights of
of a ${\mathbf g}_{-(2\A+\B)}$-submodule of $T_{r_\A r_\B}(G/B)$, so
$$\Omega (\tau_D(X,r_\A))=r_{2\A+\B} (\Omega (r_\A r_\B))=\{-(\A+\B),\A,
-(2\A+\B)\}.$$
Hence $X$ is singular at $r_\A$.
Now consider $\tau_D(X,r_\B)$ for $D=\overline{U_\A r_\A r_\B}$. By
the previous comment,
$$\Omega (\tau_D(X,r_\B))=r_\A(\Omega (r_\A
r_\B))=\{-\A,\B,-(\A+\B)\}.$$
It remains to compute $\Omega (\tau_C(X,r_\B))$ for $C=U_{\A+\B}r_\B
r_\A$.
Organizing $\Omega (r_\B r_\A)$ into $-(\A+\B)$-strings gives
$$\Omega (r_\B r_\A)=\{-\A,\B\}\cup \{\A+\B\}.$$
Since $(r_\B r_\A)^{-1}(2\A+\B)>0$, it follows from Proposition
\ref{PMG/B}
that
$$\Omega (\tau_D(X,r_\B))=r_{\A+\B}(\{-\A,-(2\A+\B),\A+\B\})=
\{-\A,\B,-(\A+\B)\}.$$
Thus, by Peterson's Theorem, $X$ is nonsingular at $r_\B$.
By the remark above, $\Omega (e)=r_\A(\Omega (r_\A))$, so
$$\Omega (e)=\{-\B,-(\A+\B),
-\A,-(2\A+\B)\}.$$
The upshot of this calculation is that the singular locus of
$X(w)$ is $X(r_\A)$.
\end{example}

\section{A Criterion For Smoothness Of $T$-varieties} \label{SC}

In this section we will prove a generalization of Theorem \ref{TE}.
Let $X$ be an irreducible $T$-variety, and let $x \in X^T$  be an
attractive $T$-fixed point. Since the action of $T$ is linearizable, and
since smoothness is a local property we may assume that $X=X_x$.
Note that we are not  assuming here that $E(X,x)$ is finite.
\begin{lemma}\label{ZNL} Let $f: X \rightarrow Y$ be a quasi-finite
equivariant morphism of $T$-varieties with $Y$ nonsingular at $f(x)$.
Let $Z \subset X$ be the ramification locus of $f$, i.e. the closed
subvariety of points, at which $f$ is not \'etale. Then either $Z$
equals
$X$, $Z$ is  empty or $Z$ has codimension one at $x$.
\end{lemma}

\begin{proof} Assume that $\codim_xZ \geq 2$. We have to show that $Z$
is
empty.  First of all, since $x$
is attractive, the image of
$f$ is contained in every connected open $T$-stable affine neighborhood
of
$f(x)$, hence in $Y_{f(x)}$. Viewing $f$ as a map to $Y_{f(x)}$, the
fibre of $f$ over $f(x)$ is finite, hence $f$ is finite. Thus, $f(X)$ is
a closed subset of the unique irreducible component of $Y_{f(x)}$
through
$f(x)$. Since $f$ is smooth somewhere it follows that $\dim X =
\dim_x Y_{f(x)}$, so $f(X)$ is the unique component of $Y_{f(x)}$
through
$f(x)$. It follows that $f(x)$ is an attractive fixed point of
$Y_{f(x)}$, and
therefore $Y_{f(x)}$ is nonsingular. Passing to  the
normalization
$\tilde X$ of $X$, we obtain an equivariant finite map $\tilde f :
\tilde X
\rightarrow Y_{f(x)}$, which is \'etale in codimension one, because the
natural map  $\tilde X \rightarrow X$ is clearly an isomorphism over $X
\setminus Z$. Thus, by the theorem of Zariski-Nagata \cite{Zar-Nag},
$\tilde f$ is
\'etale everywhere. Hence for some point $\tilde x \in \tilde X$ which
maps to $x \in X$ we have $T_{\tilde x}(\tilde X) \cong
T_{f(x)}(Y_{f(x)})$ via ${\rm d}f$. This implies that $\tilde X$ is
attractive, forcing
$\tilde f$ to be an isomorphism. Thus $f$ is birational. But being
finite, $f$ is also an isomorphism, so we are through.
\end{proof}

This gives the following criterion for smoothness of attractive
$T$-actions.

\begin{theorem}\label{GPT} Let $X$ be as above and let
$x$ be an attractive fixed point. Suppose there is a  subset $E \subset
E(X,x)$ such that every $C\in E$ is good which satisfies the following
conditions:
\begin{itemize}
\item[\rm{(1)}] \label{firstp} $|E(X,x) \setminus E| \leq \dim X - 2$.
\item[\rm{(2)}] \label{secondp} $\tau_C(X,x) = \tau_D(X,x)$ for all
$T$-curves
$C,D \in E$.
\item[\rm{(3)}] \label{fourthp} If $\tau(E)$ denotes the common value of
$\tau_C(X,x)$ for $C\in E$, then
$T_x(C) \cap \tau(E) \not = 0$ for all curves $C \in E(X,x)$.
\end{itemize}
Then $x$ is a nonsingular point of $X$.
\end{theorem}

\begin{proof} Since $x$ is attractive we may assume that
$X \subset T_x(X)$ and $x = 0$. Fix an equivariant projection $\tilde
p: T_x(X) \rightarrow \tau(E)$, and denote its restriction to $X$ by
$p$.
Since $T_x(C) \cap \tau(E) \not = 0$ for all curves
$C \in E(X,x)$, it follows that there is no $T$-curve in $p^{-1}(0)$, so
by
Lemma \ref{TC}, $\dim p^{-1}(0)=0$. This implies  $p$ is finite, since
$x$ is
attractive. Let $Z$ be the ramification locus of $p$. According to Lemma
\ref{ZNL}, we are done if $\codim_x Z \geq 2$. By assumption, if $C \in
E$, then
$C^o \subset X^*$. It follows that
$C^o \subset Z$ if and only if $dp$ has a nontrivial kernel $L
\subset T_z(X)$ for some $z \in C^o$. But then $\tau_C(L) \subset
\ker {\mathrm d}p \cap \tau(E)$. With $p$ being the projection to
$\tau(E)$,
the latter is trivial, so $\tau_C(L)$ and hence $L$ both are equal to
$0$. We
conclude that $E\cap E(Z,x)$ is empty. Thus, by condition 1), $|E(Z,x)|
\leq \dim X -2$ forcing $\dim_xZ\leq \dim X-2,$ thanks again to Lemma
\ref{TC}. This ends the proof.
\end{proof}

\begin{remark} Note that the last condition is automatically satisfied
for curves $C \in E$
since $\tau_C(C,x) \subset \tau(E)$ for such a curve. Moreover, if all
curves in $E(X,x)$ are smooth and have non collinear weights, then the
last
condition is equivalent to saying that $\tau(E) = TE(X,x)$. This in turn
implies that $E$ consists of good curves if $E \subset E(X,x)$ is a set
satisfying 2) and $|E(X,x)| =
\dim X$.
\end{remark}

We immediately conclude the following corollary, which implies the first
part of Theorem \ref{TE}.
\begin{corollary}\label{TECOR} Suppose  that either $|E(X,x)|=\dim X$ or
all
$C\in E(X,x)$ are nonsingular and any two distinct $C,D\in E(X,x)$  have
distinct tangents. Suppose also that there exist two distinct good
$T$-curves $C,D\in E(X,x)$ such that
\begin{equation}
\tau_C(X,x)=\tau_D(X,x)=TE(X,x).
\end{equation}
Then $X$ is nonsingular at $x$.
\end{corollary}

We will prove the Cohen-Macaulay  assertion in Proposition
\ref{CM}.

\begin{remark} If $X$ is normal, one does not need to assume $x$ is
attractive since
the Zariski-Nagata Theorem can be directly applied.
\end{remark}

\begin{remark} If $X$ is a $H$-variety for some algebraic group $H$, and
$(S,T)$
is an attractive slice to a $H$-orbit $Hx$ (i.e. $S \subset X$ is
locally closed, affine, stable under some
nontrivial torus $T \subset H_x$, such that $x$ is an isolated point of
$S \cap Hx$ and the natural mapping $H
\times S \rightarrow X$ is smooth at $x$), then
$T_x(Hx) \subset \tau_C(X,x)$. More precisely, one has
$\tau_C(X,x) = \tau_C(S,x)\oplus T_x(Hx)$ for all $C \in E(S,x)
\subset E(X,x)$. Thus the third condition in the theorem is always
satisfied if $E(X,x) = E(S,x) \cup E(Hx,x)$ and $E = E(S,x)$.
\end{remark}

If $X$ is a  Schubert variety $X(w)$, an explicit attractive slice for
$X$ at any $x\leq w$ is given as follows.
\begin{lemma} Let $U$ be the maximal unipotent subgroup of $B$ and $U^-$
the opposite maximal unipotent
subgroup, and suppose $x<w$. Then an attractive slice for $X(w)$ is
given
by the natural multiplication map
$$ (U\cap xU^- x^{-1})\times X(w)\cap U^- x\RA X(w).$$
\end{lemma}
We can now give a proof of Peterson's Theorem (cf.
\S 1). If $x<w$ and $\ell(w)-\ell(x)=1$, there is
nothing is to prove, since Schubert varieties are nonsingular in
codimension one.
Letting $E$ be the set of $C\in E(X,x)$ such that $C^T\subset [x,w]$,
the
existence of a slice and the
hypothesis of Peterson's Theorem imply that conditions 2) and 3) of
Theorem \ref{GPT} hold. If $\ell(w)-\ell(x)\geq
2$, then Deodhar's inequality (Proposition \ref{DI}) implies 1) holds.
Hence $X$ is nonsingular at $x$.

\section{A Fundamental Lemma}\label{FL}

In this section, $X$ will denote a $T$-variety. We will now prove a
basic
lemma which allows us to deduce good
properties of the Peterson translate from good properties of the
Peterson
translates $\tau_C(\Sigma,x)$, where
$\Sigma$ ranges over the $T$-stable surfaces containing a good $C\in
E(X,x)$.

\begin{lemma}\label{FSL} If $C \in E(X,x)$ is a good curve we have
\begin{equation*}
\tau_C(X,x)=\sum _\Sigma\tau_C(\Sigma,x)
\end{equation*}
where the sum ranges over all $T$-stable irreducible surfaces $\Sigma$
containing
$C$.
\end{lemma}
\begin{proof} Let $L \subset \tau_C(X,x)$ be a $T$-stable line. Then
by Proposition \ref{PPP} there is an $S$-line $M \subset T_z(C)$, where
$S$ is the isotropy group of an arbitrary $z
\in C^o$, such that $\tau_C(M,x) = L$. As $X$ is nonsingular at $z$,
there exists an $S$-stable curve $D$
satisfying $M\subset T_z(D)$. Setting $\Sigma = \overline{TD}$ we obtain
a $T$-stable surface, which
contains $C$, and which satisfies $L \subset \tau_C(\Sigma,x)$.
\end{proof}

Although the lemma is almost obvious, it is a great help in the case
when $X$ is a $T$-stable subvariety of $G/B$, where $G$ has no
$G_2$-factors. One reason for this is

\begin{proposition}\label{2D}
Suppose $G$ has no $G_2$-factors and let $\Sigma$ be an irreducible
$T$-stable
surface in $G/B$. Let $\sigma \in \Sigma^T$. Then $|E(\Sigma,\s)|=2$,
and
either $\Sigma$ is nonsingular at $\sigma$ or the weights of the two
$T$-curves to $\Sigma$ at $\sigma$ are orthogonal long roots $\A,\B$ in
$B_2$. In this case, $\Sigma_{\s}$ is isomorphic to a
surface of the form $z^2=xy$ where $x,y,z\in k[\Sigma_{\s}]$
 have weights $-\A,-\B,
-1/2(\A+\B)$ respectively.
In particular, if $G$ is simply laced, then
$\Sigma$ is
nonsingular.
\end{proposition}

\begin{proof} The first claim follows
easily from the fact that $\Sigma$ has a dense two dimensional
$T$-orbit (cf \cite{car-kur}).  Let $C,D$ denote the two elements
of $E(\Sigma,\sigma)$, and let $\A,\B$ denote their weights. For
any function $f \in k[\Sigma_\sigma]$ of weight $\omega$
corresponding to a $T$-line $L$ in $T_x(\Sigma)$, there is a
positive integer $N$ such that $N(-\omega) \in \Z_{\geq 0}\A +
\Z_{\geq 0} \B$. Note that the functions corresponding to $C$ and
$D$ have weights $-\A$ and $-\B$ respectively (thus, the minus
sign for $\omega$). Except for the case where $\A$, $\B$ and
$-\omega$ are contained in a copy of $B_2 \subset \Phi$ this
actually implies that $-\omega = a\A + b\B$ for suitable
nonnegative integers $a,b$. Using the multiplicity freeness of
the representation of $T$ on $k[\Sigma_\sigma]$, one is done in
these cases. In the remaining case, it turns out that $\Sigma$ is
nonsingular at $\s$ unless $\A,\B$ are orthogonal long roots in
$B_2$ (\cite{car-kur}). Let $\g=1/2(\A+\B)$. Then
$\Sigma_{\sigma}$ isomorphic to $z^2 = xy$ where $x,y,z\in
k[\Sigma_{\s}]$ correspond to $T$-lines in $T_\s(\Sigma)$ of
weights $\A,\B$ and $\g$.
\end{proof}

\section{The Span Of The Tangent Cone Of A $T$-Variety In
$G/B$}\label{tancone}

In this section and for the rest of this paper we will assume that $X$
is
a closed irreducible $T$-stable subvariety of $G/B$ and that the
underlying
semi-simple group $G$ has no
$G_2$ factors.  Recall that all
$T$-curves in $G/B$  are smooth, and two distinct $T$-curves have
different weights. Moreover, if
$C\in E(X,x)$, then $T_x(C) = {\mathcal T}_x(C) \subset {\mathcal
T}_x(X)$. In
particular, $TE(X,x)\subset \Theta_x(X)$, the $k$-linear span  of the
reduced tangent cone ${\mathcal T}_x(X)$ of $X$ at $x$.

We will now study the Peterson translates
$\tau_C(X,x)$ of $X$ along good
$T$-curves $C$ in $X$, where $x\in C^T$.  We will first show that each
$\tau_C(X,x)$ is a subspace of $\Theta_x(X)$. Hence the tangent spaces
at
$T$-fixed points behave well under
the Peterson map. On the other hand, it is an open question as to how to
explicitly describe $\Theta_x(X)$ for any
$x\in X^T$. The only relevant fact we know of is the following, proved
in
\cite{car}.

\begin{proposition} Let $S$ be an algebraic torus over $k$ and $V$ a
finite
$S$-module. Suppose $Y$ is a Zariski closed $S$-stable cone in
$V$, and let ${\cH}(Y)$ denote the convex hull of
$$\Phi(Y)=\{\A \in X(S) \mid V_\A \subset Y \}$$
in $X(S)\otimes \R$, where $V_\A$ denotes the $\A$-weight space in $V$.
Also let $\Theta(Y)$ denote the $k$-linear span of $Y$ in $V$. Then
$$\Omega (\Theta(Y))\subset {\cH}(Y).$$
\end{proposition}

In particular, if $\Phi$ is simply laced, then $\Omega
(\Theta_x(X))=\Phi \cap
{\cH}({\mathcal T}_x(X))$, so ~$\Theta_x(X)$
is  the $k$-linear span of the set of $T$-lines in ${\mathcal T}_x(X)$.
In \S \ref{SL}, we will show that if $G$ is simply laced, then, in
general,
$\Theta_x(X)=TE(X,x)$. We now prove of one of our main results.

\begin{theorem} \label{TCINTHETA}
Suppose  that $X$ is an arbitrary $T$-stable
subvariety of $G/B$. Then
for any $x\in X^T$ and any good curve $C \in E(X,x)$ we have
\begin{equation*}
\tau_C(X,x) \subset \Theta_x(X).
\end{equation*}
Moreover, if $C$ is  short, then $\tau_C(X,x)\subset  TE(X,x)$.
\end{theorem}
\begin{proof} Since $\tau_C(X,x)$ is generated by $T$-invariant surfaces
and since $\Theta_x(\Sigma) \subset \Theta_x(X)$ for all surfaces
$\Sigma
\subset X$ which contain $x$, it is enough to show the proposition when
$X$ is a surface.  If $X$ is nonsingular at $x$, then
$\tau_C(X,x)=T_x(X)$.
Otherwise we know from Proposition \ref{2D} that $X$ is a cone
over $x$, and for a cone, $\Theta_x(X) = T_x(X)$. Since $X$ is
nonsingular at
$x$ when $C$ is short, the last statement is obvious.
\end{proof}

The last assertion of the theorem gives us a generalization of
Peterson's
$ADE$ Theorem.
\begin{corollary} \label{2GC}
Let $G$ be simply laced, and suppose $X$ is rationally smooth at $x$.
Then $X$ is
nonsingular at $x$ if and only if there are two good $T$-curves in
$E(X,x)$.
Moreover, if $X$ is Cohen-Macaulay, one good $T$-curve suffices.
\end{corollary}
\begin{proof} By a result of Brion \cite{bri},  if
$X$ is rationally smooth at $x$, then $|E(X,x)|=\dim X$. By Theorem
\ref{TCINTHETA}, we have $\tau_C(X,x) =
TE(X,x)$ for every good $C\in E(X,x)$. Hence the result follows from
Theorem \ref{TE}. We will prove the last
statement below.
\end{proof}

\begin{remark}
In general, Theorem \ref{TCINTHETA} does not hold if $G=G_2$. For
example, consider
the surface $\Sigma$ given by $z^2 = xy^3$ in ${\mathbb A}^3$. Let $\A,
\B,
\g$ be characters of $T$, satisfying $\gamma = 2\A + 3\B$, and let $T$
act on
${\mathbb A}^3$ by
$$t\cdot (x,y,z)= (t^\A x,t^{(\A+2\B)}y, t^\g z).$$
Clearly $\Sigma$ is $T$-stable, and its reduced
tangent cone at $0$ is by definition $\ker {\mathrm d}z$, hence is
linear. The $T$-curve $C=\{x=0\}$ is good, and
along $C^o$, we have
$T_v(\Sigma) = \ker {\mathrm  d}x$. It follows that $\tau_C(\Sigma,x) =
\ker{\mathrm d}x$, which is
not a subspace of $\Theta_0(\Sigma)$.
It remains to remark that $\Sigma$ is open in a $T$-stable surface in
$G_2/B$, where $T$ is the usual maximal
torus and $\A$ and $\B$ are respectively the corresponding long and
short
simple roots.
\end{remark}

We now generalize Peterson's $ADE$ Theorem in a different direction.
That
is, we study which rationally smooth
$T$-fixed points of a Schubert variety (or more generally, a $T$-variety
in $G/B$) are nonsingular without the
assumption  $G$ is simply laced.  Since
$\dim_x {\mathcal T}_x(X) =\dim X$, it is  obvious that ${\mathcal
T}_x(X)$ is
linear if and only if $\dim_k \Theta_x(X) =
\dim X$. Thus, as a  consequence of Theorems \ref{TCINTHETA} and
\ref{TE}, we obtain

\begin{theorem}\label{THETAMIN}
A $T$-variety $X\subset G/B$ is nonsingular at $x \in X^T$ if and only
if
$\Theta_x(X)$ has
minimal dimension $\dim X$ and $E(X,x)$ contains at least two good
curves.
\end{theorem}
Specializing to the Schubert case gives
\begin{corollary}\label{TCL} A Schubert variety $X = X(w)$ is
nonsingular
at $x \in X^T$ if and only if all reduced tangent cones ${\mathcal
T}_y(X)$,
$x\leq y\leq w$, are linear.
Consequently, the nonsingular Schubert varieties are exactly those whose
tangent cones are linear at every $T$-fixed point.
\end{corollary}
The proof is similar to the proof of Peterson's Theorem given at the end
of \S 5 so we will omit it.

\begin{remark} Corollary \ref{TCL} gives another proof of Peterson's
$ADE$ Theorem since if
$G$ is simply laced and $X=X(w)$ is rationally smooth at $x$, then its
tangent cones are linear at every $T$-fixed
point $y$ with $x\leq y\leq w$ \cite{cp}.
\end{remark}

We now prove the second assertion of Theorem \ref{TE} that if $X$ is
Cohen-Macaulay one good $T$-curve such that $\tau_C(X,x)=TE(X,x)$
suffices to
guarantees that $x$ is nonsingular.
\begin{proposition} \label{CM} Suppose $X$ is
Cohen-Macaulay and $x \in X^T$.  $X$ is nonsingular at $x$ if and only
if
there is a good $T$-curve $C$ in $E(X,x)$ with $\tau_C(X,x) = TE(X,x)$.
\end{proposition}
\begin{proof} Let $x_1,\dots,x_n\in  k[X_x]$ be the functions
corresponding to
the $T$-curves $C_1,\dots,C_n$ in $E(X_x,x)$. Since $\dim \tau_C(X,x) =
\dim X$, and since the $T$-curves are smooth and have non collinear
weights, $n$
equals the dimension of $X$. We may assume that $C = C_1$.

Let ${\mathbf a} \subset k[X_x]$ be the ideal generated by
$x_2,\dots,x_n$. Then ${\mathbf a}$ is contained in the ideal of
$C$. Now $\tau_C(X_x,x)$ equals the span $TE(X,x)$ of the $T$-curves
$C_i$.
This means that the differentials ${\mathrm d}x_i$ are independent along
$C^o$.
As $C \cong {\mathbb A}^1$ is nonsingular, we are done if ${\mathbf
a}$ is the ideal $I(C)$ of $C$. We know that ${\mathbf a}_y =
I(C)_y$ at every stalk $k[X_x]_y$ for $y \in C^o$, since the
$({\mathrm d}x_i)_y$ are independent.

As a subset of $X_x$, $C$ is equal to the support of the Cohen-Macaulay
subscheme $Z =
\Spec(A/{\mathbf a})$ of $X_x$. Under the natural restriction, a
function $f$
on $X_x$, which vanishes on $C$,
defines a global section of ${\mathcal O}_Z$ with support contained in
$\{x\}$. It is well known that on a
Cohen-Macaulay scheme, the only such section is zero. So $f$ restricts
to
zero, and we are done.
\end{proof}

\begin{corollary} Suppose $X$ is Cohen-Macaulay and there exists a short
good $C\in E(X,x)$. Then if
$|E(X,x)|=\dim X$, then $X$ is nonsingular at
$x$. In particular, if $X$ is rationally smooth on $C$, then it is
smooth
on $C$.
\end{corollary}
\begin{proof} The only facts we have to note are, firstly, that if $C$
is
short, then $\tau_C(X,x)\subset TE(X,x)$ and, secondly, if  $X$ is
rationally smooth at $x\in X^T$, then $|E(X,x)|=\dim X$ (see
\cite{bri}).
\end{proof}

\begin{remark} Since Schubert varieties are Cohen-Macaulay, we obtain
from this yet another proof of Peterson's
$ADE$ Theorem.
\end{remark}

 For locally complete intersections or normal $T$-orbit closures in
$G/B$, Proposition \ref{CM} puts
quite strong restrictions onto the Bruhat graph $\G(X)$. In fact,
recalling that $X^*$ denotes the set of nonsingular
points of $X$, we see that if $X$ is Cohen-Macaulay,
then every rationally smooth vertex of $\G(X)$ connected to a vertex
of $\G(X^*)$ is in fact a vertex of $\G(X^*)$. Therefore we get the
following

\begin{corollary} Suppose $G$ is simply laced and $X$ is Cohen-Macaulay
, rationally smooth and
$\G(X^*)$ is
non-trivial. Then $X$ is nonsingular as long as  $\G(X)$ is connected.
\end{corollary}

\begin{remark} In fact, if $X$ is nonsingular, then its Bruhat graph is
connected.This can be shown by considering
the Bialynicki-Birula decomposition of $X$ induced by a regular element
of $Y(T)$. We will omit the details.
\end{remark}

\section{More On Peterson Translates and Schubert Varieties}
\label{PTSV}
The purpose of this section is to adress the problem that there is in
general no nice
description of $\Theta_x(X)$. Hence we would like to find a more precise
picture of $\tau_C(X,x)$ when $X$ is a Schubert variety in $G/B$, say
$X=X(w)$ and, as usual, $G$ is not allowed any  $G_2$ factors. Of
course,
if $G$ is simply laced or $C$ is short,
we've already shown $\tau_C(X,x)\subset TE(X,x)$.

It turns out that if $C$ is long, there is a $T$-subspace of
$\Theta_x(X)$,
depending only on $TE(X,x)$ and the isotropy group $B_x$ of $x$ in $B$,
which
contains most of $\tau_C(X,x)$, and the part that fails to lie in this
subspace
is easy to describe. Let ${\mathbf g}(x)$ denote the Lie-algebra of
$B_x$,
$U({\mathbf g}(x))$ its universal enveloping algebra and define
${\mathbb T}_x(X)$ to be the ${\mathbf g}(x)$-submodule
$${\mathbb T}_x(X)=U({\mathbf g}(x))TE(X,x)\subset \Theta_x(X)$$
We will show that if $C\in E(X,x)$ is both good and long and $C^T\subset
]x,w]$,
then ${\mathbb T}_x(X)$ almost contains $\tau_C(X,x)$. In fact, taking
$x=r_\A$ in Example \ref{B2} shows that in general $\tau_C(X,x)\not
\subset {\mathbb T}_x(X).$ Consequently, $\Theta_x(X)$ is not in general
equal to ${\mathbb T}_x(X).$
\begin{proposition}\label{BbbT} Assume $C\in E(X,x)$ is long, say
$C=\overline{U_{-\mu} x}$ where $\mu$ is positive and long, and suppose
$X$ is nonsingular at $y=r_\mu x$. Assume ${\mathbf g}_{\g}\subset
\tau_C(X,x)$, but
${\mathbf g}_\g
\not
\subset {\mathbb T}_x(X)$, and put $\D=\g+\mu$. Then:

\begin{itemize}

\item[\rm{(1)}] there exists a long positive root $\phi$ orthogonal to
$\mu$ such that
${\mathbf g}_{-\phi}\subset TE(X,x)$ and
$$\g=-1/2(\phi+\mu);$$

\item[\rm{(2)}]  ${\mathbf g}_{-\phi}\not \subset T_{y}(X)$;
\item[\rm{(3)}] if $\D>0$, then $x^{-1}(\D)<0$ and
$${\mathbf g}_\g \oplus {\mathbf g}_{\D}\oplus {\mathbf g}_{-\mu}
\subset \tau_C(X,x),\quad  {\mathbf g}_\g \oplus {\mathbf
g}_{\D}\oplus {\mathbf g}_\mu \subset T_{y}(X);$$

\item[\rm{(4)}] on the other hand, if $\D<0$, then $x^{-1}(\D)>0$ and
$${\mathbf g}_\g \oplus {\mathbf g}_{-\D}\oplus {\mathbf g}_{-\mu}
\subset \tau_C(X,x), \quad {\mathbf g}_{-\g} \oplus {\mathbf
g}_{\D}\oplus {\mathbf g}_\mu \subset T_{y}(X);$$

\item[\rm{(5)}] In particular, if $D=\overline{U_{-\phi}x}$, then
$\tau_C(X,x)\neq\tau_D(X,x)$.

\end{itemize}

\end{proposition}
\begin{proof}
The existence of a long root $\phi$ satisfying all the conditions
 in part (1) except possibly positivity
follows from the Fundamental Lemma (\ref{FSL}), Proposition
\ref{2D} and the fact that ${\mathbf g}_{-\mu} \subset TE(X,x)$.
To see $\phi$ is positive suppose otherwise. It's then clear that
$\D\in \Phi^+$. If
$x^{-1}(\D)>0$ also, then ${\mathbf g}_\g \subset {\mathbb T}_x(X)$
since
$\g=-\mu+\D$, contradicting the assumption. Hence $x^{-1}(\D)<0$. But as
$\tau_C(X,x)$ is a ${\mathbf g}_\mu$-module, it follows immediately that
${\mathbf g}_\g \oplus {\mathbf g}_\D \subset \tau_C(X,x)$. Since $\mu$
is long, Proposition
\ref{PMG/B} implies
$${\mathbf g}_\g\oplus {\mathbf g}_{\D}\subset T_y(X).$$
Moreover, since $x^{-1}(\phi)=y^{-1}(\phi)>0$, it also follows
that $ {\mathbf g}_{-\phi} \subset T_y(X)$, so $\mu,~\D,~-\phi$
constitute
a complete $\g$-string  occuring in $\Omega (T_y(X))$. Since $X$ is
nonsingular at $y$ and $\g,y^{-1}(\g)<0$,
we get the inequality
$y<r_\g y\leq w$. Thus
 $X$ is nonsingular
at $r_\g y$. Letting $E$ be the
good $T$-curve in $X$ such that $E^T=\{y,r_\g y \}$, we have
$\tau_E(X,y)=T_y(X)$, so the string $\mu,~\D,~-\phi$ also has to occur
in
$\Omega (T_{r_\g y}(X))$. In particular, ${\mathbf g}_{-\phi}
\subset  TE(X, r_\g y)=T_{r_\g y}(X)$, and hence $r_\phi r_\g
y\leq  w$. But this means
$$r_\g x=r_\g r_\mu y=r_\g r_\mu r_\g r_\g y=r_\phi r_\g y\leq w,$$
so ${\mathbf g}_\g \subset TE(X,x)$.  This is a contradiction, so
$\phi>0$.

We next prove (2).  Recall $y=r_\mu x$ and suppose to the contrary that
${\mathbf g}_{-\phi} \subset
T_{r_\mu x}(X)$. If $\D>0$, we can argue exactly
as above, so we are reduced to assuming $\D<0$. If
$x^{-1}(\D)<0$, then ${\mathbf g}_\g \subset {\mathbb T}_x(X)$ due
to the fact that $\g=-(\phi+\D)$. Thus $x^{-1}(\D)>0$, whence
$y^{-1}(\g)>0$ and
${\mathbf g}_{-\g}\subset T_y(X)$ since
$-\g>0$. Moreover, ${\mathbf g}_\D\subset T_y(X)$ since ${\mathbf g}_\g
\subset
\tau_C(X,x)$. Hence, we have
$${\mathbf g}_\mu \oplus {\mathbf g}_{-\phi} \oplus {\mathbf g}_\D
\oplus {\mathbf
g}_{-\g}
\subset T_y(X).$$
Now $\D$ and $-\g$ form a $\phi$-string in $\Omega (T_y(X))$, and since
$y<r_\phi y
\leq w$, it follows as above that ${\mathbf g}_\D \subset T_{r_\phi
y}(X)$.
Therefore, $r_\D r_\phi y\leq w$. But $r_\g x=r_\D r_\phi y$, so we have
a
contradiction. Hence, ${\mathbf g}_{-\phi}\not \subset T_{r_\mu x}(X)$.
To prove (3), note that, as usual, $x^{-1}(\D)<0$, so ${\mathbf g}_\D
\subset
\tau_C(X,x)$ by the ${\mathbf g}_\mu$-module property. To obtain (4),
note
that if
$\D<0$, then $x^{-1}(\D)>0$, so $y^{-1}(-\g)<0$. As $\g<0$, this implies
that
${\mathbf g}_{-\g}\subset T_y(X)$, so in fact
${\mathbf g}_{-\g} \oplus {\mathbf g}_\D \oplus {\mathbf g}_\mu \subset
T_y(X)$. The proof of (5) is clear, so the proof is now finished.
\end{proof}

Now fix $C$ and $\mu\in \Phi^+$ as above and let $I_\mu \subset \Phi$
consist of all negative $\gamma$ such that:
\begin{itemize}
\item[(a)] $\gamma=-1/2(\mu + \phi)$, where $\phi$ satisfies conditions
(1) and (2) of Proposition \ref{BbbT},

\item[(b)] $\delta =\mu +\gamma\in \Phi$, and
\item[(c)] $\delta$ satisfies conditions (4) and (5).
\end{itemize}
Put $V_C=\bigoplus_{\g\in I_\mu}{\mathbf g}_\g$. Notice that
$V_C\subset T_x(G/B)$. Proposition \ref{BbbT} thus gives  the following:

\begin{corollary} Assuming the previous hypotheses, we have
$$\tau_C(X,x)\subset {\mathbb T}_x(X)+V_C\subset \Theta_x(X).$$
\end{corollary}

\section{The Simply Laced Case}\label{SL}

The crucial point in the proof of Peterson's $ADE$ Theorem is the fact
that
every $T$-stable line in the span of the tangent cone of $X$ comes
from a $T$-curve in $X$. It turns out that this is true for any closed
$T$-variety $X\subset G/B$ as long as $G$ is simply laced. We now prove
this fact.
\begin{theorem} \label{TC=TL} Suppose $G$ has no $G_2$-factors. Let $L
\subset \Theta_x(X)$ be a $T$-stable line with weight $\omega$. Then
\begin{equation*}
\omega = \frac{1}{2}(\A + \B)
\end{equation*}
where $\A$ and $\B$ are the weights of
suitable $T$-curves $C$ and $D$, respectively.
If, moreover, $G$ is simply laced, then $\A = \B = \omega$,
hence $L$ is the tangent line of a $T$-curve $C \in E(X,x)$.
\end{theorem}

\begin{proof}
We will prove the following equivalent 'dual' statement: if $\omega$ is
the weight
of a function corresponding to a $T$-stable line $L \subset
\Theta_x(X)$, then there are $\A$ and $\B$
with $\omega = 1/2(\A + \B)$ where $\A$ and $\B$ are the weights of
functions corresponding to
$T$-curves $C$ and $D$, respectively.

Let $z \in k[X_x]$ be the $T$-eigenfunction corresponding to $L$, and
let
$$x_1,x_2,\dots,x_n \in k[X_x]$$
be those corresponding to the $T$-curves
$C_1,C_2, \dots, C_n$ through $x$. Consider the unique linear
projections
$$\tilde x_i:T_x(X) \RA T_x(C_i),\quad
\tilde z:T_x(X) \RA L$$
which restrict respectively to $x_i, z \in k[X_x]$.

Since the (restriction of the) projection $X_x \RA \bigoplus
T_x(C)$ has a finite fibre over $0$, $k[X_x]$ is a finite
$k[x_1,x_2, \dots , x_n]$-module. In particular $z \in k[X_x]$ is
integral over $k[x_1,\dots,x_n]$. We obtain a relation
\begin{equation} \label{INTEQ}
z^ N = p_{N-1} z^{N-1} + p_{N-2}z^{N-2} + \dots + p_1 z
+ p_0
\end{equation}
where $N$ is a suitable integer and the $p_i\in k[x_1,\dots,x_n]$.
Without loss of generality we may assume that
every summand on the right hand side is a $T$-eigenvector with
weight $N \omega$. Let $P_i \in k[\tilde x_1, \dots, \tilde x_n]$ be
polynomials restricting to $p_i$, having the same weight $(N-i)\omega$
as $p_i$. Then every monomial $m$ of $P_i$ has this weight too.
If for all $i$ every such monomial $m$ has degree $\deg m > N - i$, then
$p_i z^{N- i}$ is an element of ${\mathbf m}^{N + 1}_x$, where
${\mathbf m}_x$ is the maximal ideal of $x$ in $k[X_x]$. This means that
$\tilde z$ vanishes on the tangent
cone of $X_x$, so $L \not \subset \Theta_x(X)$,
which is a contradiction.

Thus, there is an $i$ and a monomial $m$ of $P_i$, such that $\deg
m \leq M = N - i$. Let $m = c \tilde x_1^{d_1} \tilde x_2 ^{d_2}
\dots \tilde x_n^{d_n}$, with integers $d_j$ and a nonzero $c \in
k$. So $\sum_j d_j \leq N$. Let $\A_j$ be the weight of $\tilde
x_j$. Then we have
\begin{equation*}
 M \omega = \sum d_j \A_j
\end{equation*}
After choosing a new index, if necessary, we may assume that $d_j
\not = 0$ for all $j$. Let $F$ be a nondegenerate bilinear form on
$X(T)\otimes {\mathbb Q}$ which induces the length function on
$\Phi$. We have to consider two cases. First suppose that $\omega$ is
a long root, with length say $L$. Then $F(\A_j, \omega) \leq L^2$ with
equality if and only if $\A_j = \omega$. Thus, $M L^2 = \sum d_j
F(\A_j,\omega) \leq M \max_{j} F(\A_j, \omega) \leq M L^2$ and so there
is
a $j$ with $\A_j = \omega$ and we are done, since this implies $\tilde
z = \tilde x_j$. Note that, although we are considering all roots short
in case $G$
is simply laced, this contains actually the case that all roots have the
same length.

Now suppose $\omega$ is short, having length $l$. In this case
$F(\A_j, \omega) \leq l^2$. Since $M l^2 =  MF(\omega,\omega) = \sum_j
d_j
F(\A_j, \omega)$ and since $\sum d_j \leq M$, it follows that all
$\A_j$ satisfy $F(\A_j, \omega) = l^2$. If there is a $j$ such that
$\A_j = \omega$, then, as above, we are done. Otherwise for each $j$,
$\A_j$ is long, and $\A_j$ and $\omega$ are contained in a copy $B(j)
\subset \Phi$ of $B_2$. There is a long root $\B_j \in B(j)$ with
$\A_j + \B_j = 2\omega$. We have to show that there are $j_0$ and
$j_1$ so that $\B_{j_0} = \A_{j_1}$. Fix $j_0 = 1$ and let $\A =
\A_1$, $\B = \B_1$. Then $F(\A,\B) = 0$.
 This gives us the result: $M l^2 = MF(\omega, \B) = 0 + \sum_{j>1}
F(\A_j, \B)$. Now if all $F(\A_j,\B)$ are less or equal $l^2$, this
last equation cannot hold, since $\sum_{j>1}d_j
< M$. We conclude that there is a $j_1$ so that $F(\A_{j_1}, \B) =
L^2$, hence $\A_{j_1} = \B$, and we are through.
The statement for $G$ simply laced follows from this, since
in this case all roots have the same length, hence are long.
\end{proof}

For completeness, we state the following corollary.
\begin{corollary}\label{}Suppose the $G$ is simply laced and that $X$
is a $T$-variety in $G/B$ such that $\dim X\ge 2$. Then $X$ is smooth at
$x\in X^T$ if and only if
$|E(X,x)|=\dim X$ and there are at least two good $T$-curves at $x$.
\end{corollary}

Theorem \ref{TC=TL} implies in particular that the linear spans of the
tangent cones of two
$T$-varieties behave nicely under intersections, and this allows us to
deduce a somewhat surprising
fact about the intersection of the tangent spaces of two
$T$-varieties at a common nonsingular point.
\begin{corollary} \label{INT} Suppose the $G$ is simply laced and that
$X$
and $Y$ are $T$-varietes in $G/B$. Suppose also that $x\in
X^T\cap Y^T$. Then
$$\Theta_x(X\cap Y)=\Theta_x(X)\cap \Theta_x(Y).$$ Furthermore, if both
$X$ and $Y$ are
nonsingular at $x$, then
$$T_x(X\cap Y)=T_x(X)\cap T_x(Y).$$
In particular, if $|E(X\cap Y,x)|=\dim X\cap Y$, then $X\cap Y$ is
nonsingular at
$x$.
\end{corollary}

\begin{proof} The first assertion is a consequence of Theorem
\ref{TC=TL} and the fact that $E(X,x)\cap
E(Y,x)=E(X\cap Y, x)$. For the second, use the fact that if a variety
$Z$ is nonsingular at $z$, then
$T_z(Z)=\Theta_z(Z)$. Thus
\begin{eqnarray*}
T_x(X)\cap T_x(Y)&=&\Theta_x(X)\cap \Theta_x(Y)\\
                  &=&\Theta_x(X\cap Y)\\
                  &\subset & T_x(X\cap Y)\\
                  &\subset & T_x(X)\cap T_x(Y)
\end{eqnarray*}
The final assertion follows from the fact that if $X$ and $Y$ are both
nonsingular at $x$, then $T_x(X)\cap T_x(Y)=TE(X\cap Y,x)$.
\end{proof}

For example, it follows that the in the simply laced setting, the
intersection of a Schubert variety
$X(w)$ and
a dual Schubert variety $Y(v)=\overline{B^-v}$ is nonsingular at each
$x$ with $v\leq
x \leq
w$ as long as each of the constituents is nonsingular at $x$.

\begin{remark} The previous corollary was stated in \cite{cp} (cf.
Theorem H and
Corollary H) for so called shifted Schubert varieties, that is any
subvariety of $G/B$ of the form
$X(y,w)=yX(w)$, where $y,w\in W$. One can in fact say a little more for
shifted
Schubert varieties in type
$A$ because from a result of Lakshmibai and Seshadri \cite{ls} saying
that if $G$ is of type $A$, then
$\Theta_x(X)=T_x(X)$ for every shifted Schubert variety $X$. Namely,
$$T_x(X)\cap T_x(Y)=T_x(X\cap Y)$$
for any two shifted Schubert varieties $X$ and $Y$ meeting at $x\in W$.
\end{remark}

\section{Singular Loci of Schubert Varieties}\label{SLSV}

In this section we give an algorithm for
computing the singular locus $X^{\times}$  of a Schubert variety
$X=X(w)$
assuming $G$ has no $G_2$ factors. Obviously $X^{\times}$ is a union of
Schubert varieties, so we only have to compute the maximal elements of
$X^{\times T}$. Schubert varieties $X=X(w)$ being Cohen-Macaulay , we
know $x<w$ is a nonsingular point as long as
$E(X,x)$ contains a short good $T$-curve and  $|E(X,x)|=\dim X$. Hence
maximal elements $x$ of $X^{\times}$
either have the property that  $|E(X,x)|>\dim X$ or every good $T$-curve
at $x$ is long.

On the other hand, we can use Proposition \ref{PMG/B}, which gives a
criterion for
deciding when $\tau_C(X,x)=\tau_D(X,x)$ that doesn't depend on knowing
either or both
of $C,D\in E(X,x)$ are good. Suppose $C=\overline{U_\A x}$ and
$D=\overline{U_\B z}$ where $\A,\B>0$. Let $y=r_\A
x, z=r_\B x$ so $y,z\in [x,w]$. By Proposition \ref{PMG/B},
$\tau_C(X,x)=\tau_D(X,x)$ if and only if
$r_\A\Omega (T_{y}(X)^*)=r_\B\Omega (T_{z}(X)^*)$,
or,  equivalently,
$\Omega (T_{y}(X)^*)=r_\A r_\B\Omega (T_{z}(X)^*)$, where $T_{y}(X)^*$
is the ${\mathbf
g}_{-\A}$-submodule of $T_y(G/B)$ defined in \S \ref{G/B}. Note that we
are
working in $T_y(G/B)$ instead of $T_x(X)$. Now
${\mathrm d}\dot{r_\A}{\mathrm d}\dot{r_\B}(T_{z}(X)^*)$ is a ${\mathbf
g}_{r_\A(\B)}$-module, and consequently, this implies that
$T_{y}(X)^*$ is a module for the subalgebra of ${\mathbf g}$ generated
by
${\mathbf g}_{-\A}$ and ${\mathbf g}_{r_\A(\B)}$.

\begin{proposition}\label{C=D} Assuming the notation is as above,
$\tau_C(X,x)=\tau_D(X,x)$ if and only if $T_y(X)^*$ is a ${\mathbf
g}_{r_\A(\B)}$-submodule of  $T_y(G/B)$ whose leading weights are  the
$r_{\A}r_{\B}(\gamma)$, where
$\gamma$ runs through the set of  leading weights for the
${\mathbf g}_{-\B}$-module $T_z(X)^*$. Moreover, if
$C_1, C_2, \dots C_k\in E(X,x)$ where each $C_i^T=
\{x,y_i\}$ with  $y_i=t_ix>x$ and if all $\tau_{C_i}(X,x)$ coincide,
then every $T_{y_i}(X)^*$ is  a
module for the subalgebra ${\mathbf m}_i$ of ${\mathbf g}$
generated by
$${\mathbf g}_{t_i(\A_1)}\oplus {\mathbf g}_{t_i(\A_2)}\oplus \cdots
\oplus
{\mathbf g}_{t_i(\A_{t_k})}.$$
\end{proposition}
 The last assertion is a consequence of the Jacoby identity. Note that
in this result, there is no assumption that the $T$-curves be good.

The algorithm for determing $X^{\times}$ is now clear. Suppose one knows
that $y\in X^{*T}$, and assume $x=r_\A y<y$ where $\A >0$. Clearly if

$${\mathrm d}\dot{r_\A}(\Omega (T_y(X)^*))\neq \{\g \mid
x^{-1}(\g)<0,r_\g x\leq
w\},$$
then $x\in X^{\times}$. If equality holds, then it suffices to apply
Proposition \ref{C=D} to any good $D\in E(X,x)$.
Thus the algorithm requires checking whether $X$ is nonsingular at any
$z\in X^{*T}$ with $z>x, z\neq y$ and $sz=x$ for some $s\in R$.

\section{Generalizations to $G/P$}\label{G/P}
As usual, assume $G$ is semi-simple and has no $G_2$ factors, and
suppose
$P$ is a parabolic subgroup of $G$ containing
$B$. In this section, we will indicate which results extend to
$T$-varieties in
$G/P$. Let $\pi:G/B\RA G/P$ be the natural projection.
The extensions to $G/P$ are based on the following lemma.
\begin{lemma}\label{GmodP} Let $Y\subset G/P$ be closed and $T$-stable,
and put $X=\pi^{-1}(Y)$.
Then:
\begin{itemize}
\item[\rm{(1)}] the projection $\pi:X\RA Y$ is a smooth morphism, hence
$X^*=\pi^{-1}(Y^*)$;

\item[\rm{(2)}] for all $x\in X^T$, ${\rm d}\pi_x :T_x(X)\RA T_y(Y)$ is
surjective and
$${\rm d}\pi_x(\Theta_x(X))=\Theta_y(Y),$$
where $y=\pi(x)$;

\item[\rm{(3)}] $\pi(E(X,x))=E(Y,y)$; and

\item[\rm{(4)}] if $C\in E(X,x)$ is good and $\pi(C)$ is a curve, then
$\pi(C)\in E(Y,y)$ is good and
$${\rm d}\pi_x(\tau_C(X,x))=\tau_{\pi(C)}(Y,y).$$
\end{itemize}
\end{lemma}

\begin{proof} The first statement (1) is standard.
Moreover, ${\rm d}\pi_x$ is surjective for  all $x\in X$ and
${\rm d}\pi_x$ maps the schematic tangent
cone of $X$ at $x$ onto that of $Y$ at $y$.
Consequently, it is also
a surjection of the  the associated reduced varieties.
Since ${\rm d}\pi_x$ is linear, (2) is established. (3) is
an immediate consequence of the fact that
$\pi(E(G/B,x))=E(G/P,y)$. Part (4) follows from the
existence of a local $T$-equivariant
section of $\pi$, the smoothness of $\pi$ and Lemma \ref{FSL}.
\end{proof}

\begin{corollary}\label{NN} Assume $G$ has no $G_2$
factors, and suppose $Y$ is any $T$-variety in $G/P$.
If $\dim Y \ge 2$, then
$Y$ is smooth at the $T$-fixed point $y$ if and only if
$E(Y,y)$ contains two good $T$-curves and the
reduced tangent cone to $Y$ at $y$ is linear.
\end{corollary}

\begin{proof}
Apply the previous lemma and Theorem \ref{THETAMIN} to $X=\pi^{-1}(Y)$
at
any $x\in \pi^{-1}(y)^T$, which,
by the Borel Fixed Point Theorem, is
non-empty since $y\in Y^T$.
\end{proof}

If $G$ is simply laced, there is more.
\begin{corollary} Assume $G$ is simply laced. Then for any
$T$-variety $Y$ in $G/P$ and $y\in Y^T$,
$$\Theta_y(Y)=TE(Y,y).$$
In particular, if $\dim Y\ge 2$, then $Y$ is smooth at $y$ if and only
if $|E(Y,y)|=\dim Y$ and $y$ lies on at
least two good $T$-curves.
\end{corollary}

We also have
\begin{theorem}\label{RSG/P} If $G$ is simply laced, then every
rationally
smooth $T$-fixed point of a Schubert variety $Y$ in
$G/P$ is nonsingular.
\end{theorem}
\begin{proof} Let $y\in Y^T$ be a  rationally
smooth $T$-fixed point of $Y$. Using the relative order,
we may without loss of generality assume
that if $z\in Y^T$ and $z>y$, then $Y$ is nonsingular
at $z$. By the relative version of Deodhar's Inequality
and the fact that the singular locus
of $Y$ has codimension at least two (as $Y$ is normal), there are at
least
two good $T$-curves in
$E(Y,y)$.  Since $|E(Y,y)|=\dim Y$ (\cite{bri}), the proof is done.
\end{proof}

Finally, we state a $G/P$ analog of Corollary \ref{RSGmodP}.
\begin{corollary} If $G$ is simply laced, a Schubert variety in $G/P$ is
nonsingular
if and only if
the  Poincar\'e polynomial of $Y$ is symmetric
if and only if $|E(Y,y)|=\dim Y$ for every $y\in Y^T$.
\end{corollary}

\pagebreak
\section{ A Remark and Two Problems} \label{QAR}

Although we have not yet given an explicit example, it is definitely not
true that
 in the simply laced setting, every rationally smooth $T$-variety in
$G/B$ is nonsingular. In fact,
there are $T$-orbit closures in types $D_n$ if
$n>4$ and in $E_6, E_7, E_8$ which are rationally smooth but non-normal,
hence singular. For more information, see \cite{car-kur,mo}.
A final comment is that one of the most basic open problems about
Schubert varieties in our context is to describe the $T$-lines in the
linear span of the tangent
cone at a $T$-fixed point in the non-simply laced setting.
Once this is
settled, we will have a complete picture of the singular loci of all
Schubert varietes.
Another unsolved problem is to identify all the $T$-lines in the tangent
space
at at $T$-fixed point. There are results in this
direction in papers of Lakshmibai and Seshadri \cite{ls}, Lakshmibai
\cite{la} and Polo \cite{po}. The
natural  conjecture that these tangent spaces are spanned by Peterson
translates is, in
light of Theorem 1.3, seems to  be  only true for type $A$.

\newpage
{\footnotesize
\begin{theckbibliography}{00}

\bibitem{bri} M\ Brion: {\it Rational smoothness and fixed points of
torus actions}
Transformation Groups, Vol.\ 4, No.\ 2-3, (1999) 127-156

\bibitem{cp}  J B\ Carrell: {\em The Bruhat Graph of a Coxeter Group, a
Conjecture of Deodhar, and Rational Smoothness of Schubert Varieties.}
Proc.\ Symp.\ in Pure Math.\ {\
bf 56}, No.\ 2, (1994), Part 1, 53-61.

\bibitem{spsv} J B\ Carrell: {\em On the smooth points of a
Schubert variety.}
Representations of Groups, Canadian
Mathematical Conference Proceedings {\bf 16}
(1995), 15-35.

\bibitem{car} J B\ Carrell: {\em The span of the tangent cone of a
Schubert
variety.} Algebraic Groups and Lie Groups, Australian Math. Soc. Lecture
Series
{\bf 9}, Cambridge Univ. Press (1997), 51-60.

\bibitem{car-kur} J B\ Carrell and A\ Kurth: {\em Normality of Torus
Orbit
Closures on $G/P$}, preprint (1998).

1-25.

\bibitem{deod} V\  Deodhar: {\it Local Poincar\'e duality and
nonsingularity of Schubert
varieties}, Comm. in Algebra {\bf 13} (1985), 1379-1388.

\bibitem{Zar-Nag}A\ Grothendieck: {\it Cohomologie Locale des Faisceaux
Coh\'erents et Th\'eor\`emes de Lefschetz Locaux et Globaux}, SGA {\bf
2} (1962), X, 3.4

\bibitem{kl1} D\ Kazhdan and G\ Lusztig:  {\em Representations of
Coxeter groups and Hecke algebras}, Invent. Math. {\bf 53} (1979),
165-184.

\bibitem{kl2}  D\ Kazhdan and G\ Lusztig:  {\em Schubert varieties and
Poincar\'e duality}, Proc. Symp. Pure Math. A.M.S. {\bf 36} (1980),
185-203.


\bibitem{kut} J\ Kuttler: {\it Ein Glattheitskriterium F$\ddot{u}$r
Attraktive Fixpunkte von Torusoperationen}, Diplomarbeit, U. of Basel
(1998).

\bibitem{la} V\ Lakshmibai: {\it Tangent spaces to Schubert varieties},
Math.
Res. Lett. 2(1995), no. 4, 473--477.

\bibitem{ls} V\ Lakshmibai and C\ Seshadri:  {\it Singular locus of a
Schubert variety}, Bull. Amer. Math. Soc. (N.S.) {\bf 5} (1984),
483-493.

\bibitem{mo} J\ Morand: {\it Closures of torus orbits in adjoint
representations of
semisimple groups}, C.R.Acad. Sci. Paris Se\'r. I Math. {\bf 328}(1999),
no. 3, 197-202.

\bibitem{mu} D\ Mumford: {\em The Red Book of Varieties and Schemes.}
Springer LN {\bf 1358}, Springer-Verlag, Berlin Heidelberg (1988).

\bibitem{no}  A\ Nobile: {\em Some properties of the Nash blowing up},
Pac. Jour. Math. {\bf 60} (1975), 297-305.



\bibitem{po} P\ Polo:  {\it On Zariski tangent spaces of Schubert
varieties and a proof of a conjecture of Deodhar}, Indag. Math. {\bf 11}
(1994), 483-493.

\bibitem{sumi1} H\ Sumihiro: {\it Equivariant Completion}, J. Math.
Kyoto Univ. {\bf 14} (1974), 1-28
\bibitem{sumi2} H\ Sumihiro: {\it Equivariant Completion II}, J. Math.
Kyoto Univ. {\bf 15} (1975), 573-605

\end{theckbibliography} }

%

\pagebreak

{\scriptsize
\begin{tabular}{l}
James B.\ Carrell\\
Department of Mathematics\\
University of British Columbia\\
Vancouver, Canada V6T 1Z2\\
E-mail.\ carrell@math.ubc.c{a}\\\\
Jochen Kuttler\\
Department of Mathematics\\
University of California at San Diego\\
LaJolla, CA 92093\\
and\\
Mathematisches Institut\\ 
Universit\"at Basel\\
Rheinsprung 21\\
CH-4051 Basel\\
Switzerland
E-mail.\ kuttler@math.unibas.ch\\
\end{tabular}}

\end{document}